\documentclass[12pt]{article}
\usepackage{amssymb,latexsym,amsfonts,amsmath,graphicx}
\newcommand{\const}{{\rm const}}
\newcommand{\mynorm}{\mathcal{M}}

\begin{document}
\title{Implementation for blow up of tornado-type solutions for complex
version of $3D$ Navier-Stokes system.}
\author
{Arnold M.D. \thanks{International institute of the Earthquake Prediction Theory, Russian Academy of Sciences.} \thanks{Institute of the Information Transition Problems, Russian Academy of Sciences}
\and
Khokhlov A.V.\footnotemark[1]
}
\maketitle

\abstract{\noindent         
We consider Cauchy problem for Fourier transformation of $3$-dimensional Navier-Stokes system with zero external force. Using initial data purposed by Dong Li and Ya.G.Sinai in \cite{DS} we implement self-similar regime producing fast growing behavior of the energy of solution while time tends to critical value $T_{cr}$.}

\section {Theoretical results used.}
Fourier transform of the $3D$ Navier-Stokes system with purely imaginary initial data can be written in the form 
\begin{equation}
\label{eq: NSS}
v(k,t)=e^{-t|k|^2}v(k,0)+\int\limits_0^te^{-(t-s)|k|^2}\int\limits_{\mathbb{R}^3}\langle
k,v(k-l,s)\rangle P_k(v(l,s))dlds
\end{equation}
where $k\in\mathbb{R}^3$ corresponds to Fourier mode of the solution, $t\in [0, T_{cr}]\subset \mathbb {R}_+$ is time, $v(k,t)$ corresponds to imaginary part of the solution,
$\langle \cdot,\cdot \rangle$ denotes Euclidian inner product. $P_k$ is a Leray projection on the subspace orthogonal to $k$ and has the form 
$P_k=Id-\dfrac{\langle k,\cdot\rangle}{|k|^2}k$. Viscosity supposed to be equal $1$, external forcing is put to zero. Incompressibility condition means that for all $t\in[0,T_{cr}]$
\begin{equation} \label{eq: incompress} 
v(k,t)\bot k \qquad \mbox{for any} \, k\in \mathbb{R}^3
\end{equation}

For given function $c_0(k)$ consider one-parametric family of initial conditions  
$v_A(k,0)=Ac_0(k)$.
Accordingly to Ya.G. Sinai \cite{S1}, \cite{S2} for any $c_0$ with fast decay at infinity, for every $A$ there exists such constant $T_{cr}$ that for $t\in[0,T_{cr}]$ there exists unique solution of the system \eqref{eq: NSS} which can be represented as a power series of the parameter $A$:
\begin{equation}
v_A(k,t)=Ae^{-t|k|^2}c_0(k)+\int\limits_0^te^{-(t-s)|k|^2}\sum\limits_{p>1}A^p
h_p(k,s)ds\label{eq: expansion} 
\end{equation}
 
Substituting \eqref{eq: expansion} into \eqref{eq: NSS} one can express coefficient $h_p(k,s)$ trough coefficients with fewer indices and from initial data $c_0$ with the following recurrent relations:

\begin{gather}
h_2(k,s)=\int\limits_{\mathbb{R}^3}\langle l, c_0(k-l)\rangle P_k
c_0(l) \exp\{-s(|l|^2 + |k-l|^2 \}dl\notag\\
h_p(k,s)=\int\limits_0^s ds_1\int\limits_{\mathbb{R}^3}\langle l, c_0(k-l)\rangle P_k
h_{p-1}(l,s_1) \exp\{-(s-s_1)|l|^2 - s|k-l|^2 \}dl+\notag\\
+\int\limits_0^s ds_2\int\limits_{\mathbb{R}^3}\langle l, h_{p-1}(k-l,s_2)\rangle P_k
c_0(l) \exp\{-s|l|^2 - (s-s_2)|k-l|^2 \}dl+\notag\\
+\sum\limits_{p_1+p_2=p} \int\limits_{0}^s ds_1 \int\limits_0^s
ds_2 \int \limits_{\mathbb{R}^3} \langle l,
h_{p_1}(k-l,s_1)\rangle P_k h_{p_2}(l,s_2)
\\\notag\exp\{-|k-l|^2(s-s_1)-|l|^2(s-s_2)\}dl
\end{gather}

Consider initial data with compact support $C$ located in the ball of radius $r$ with the center at the point $(0,0,R)$. From recurrent relations (4) it follows that the coefficient $h_p$ is supported in the algebraic sum $\underbrace{C+C+\cdots +C}\limits_{p\, \mathrm{times}}$. Then one can fix $1<r$ such that for the ratio $\frac R r$ large enough supports of the coefficients $h_p$ and $h_{p'}$ $p'\ne p$ does not intersect. In more details, support of each function $h_p$ lies in the ball of radius $O(1)\sqrt {p}$ with center at the point $(0,0, R\cdot p)$.    
 
We reproduce sketch of the arguments of Ya.G. Sinai and Dong Li \cite{DS}. Third component of coefficient $h_p$ can be expressed through first and second ones due to the incompressibility condition \eqref{eq: incompress}. Since all terms in (4) where either $p_1$ or $p_2$ is small comparable to $\sqrt{p}$ can be estimated with some function which decays with $p$ faster than exponent one can concider only the part where both $p_1\sim \sqrt{p}$ and $p_2\sim \sqrt{p}$. Change variables according to the ball containing the support of function $h_p$ and consider a function 
$g_{p}(y,s)=h_p((0,0,R)+\sqrt{R}y,s)$. After renormalisation, we get  
\begin{equation}
\label{eq: renormalize}
\begin{split}
&g_{p}(y,s)=\frac {p^{\frac
52}}{R^2}\sum\limits_{p_1+p_2=p\atop{p_1,p_2>\sqrt{p}}}\frac{p^2}{p_1^2p_2^2}
\int\limits_{\mathbb{R}^3}g_{p_2}(x\sqrt{\frac{p}{p_2}},s)\left(\frac{p_2}{p_1}(y_1-x_1)g_{p_1}^{(1)}((y-x)\sqrt{\frac{p}{p_1}},s)+\right.
\\
 &\left.+
(y_2-x_2)g_{p_1}^{(2)}((y-x)\sqrt{\frac{p}{p_1}},s)+x_1
g_{p_1}^{(1)}((y-x)\sqrt{\frac{p}{p_1}},s)+x_2g_{p_1}^{(2)}((y-x)\sqrt{\frac{p}{p_1}},s)\right)dx
\end{split}
\end{equation}

Write down functions $g_p(y,s)$ in the form  
\begin{equation}
\label{eq: gp}
\begin{split}
g_p(y,s)=\Lambda^p(s)p \sigma_1
\exp\left\{-\frac{\sigma_1}{2}(|y_1|^2+|y_2|^2)\right\}\sqrt{\sigma_2}
\exp\left\{\sigma_2|y_3|^2\right\}\\
\left(H_1(y)+\delta_1^{(p)}(y,s),
H_2(y)+\delta_2^{(p)}(y,s),\delta_3^{(p)}(y,s)\right)
\end{split}
\end{equation}

While $p$ tends to infinity one can show that functions $\delta_1^{(p)}(y,s)$, $\delta_2^{(p)}(y,s)$ and $\delta_3^{(p)}(y,s)$ became $o(1)$.
Sum in the right-hand side of the \eqref{eq: renormalize} can be considered as integral sum with variable  $\gamma=\dfrac{p_1}{p}$. For the components $H_1$ and $H_2$ we get from \eqref{eq: gp},\eqref{eq: renormalize} (for details see \cite[Eq.21]{DS}))
\begin{equation}
\label{transformedeq}
\begin{split}
(H_1(x), H_2(x))= \frac{2\pi}{\sigma^{(1)}}&\exp\left\{\frac{\sigma^{(1)}|x|^2}{2}\right\}\int\limits_0^1d\gamma
\int\limits_{\mathbb{R}^2}
\frac{(\sigma^{(1)})^2}{(2\pi)^2\gamma(1-\gamma)}\exp\left\{-\frac{\sigma^{(1)}|x-y|^2}{2\gamma}-
\frac{\sigma^{(1)}|y|^2}{2(1-\gamma)} \right\}\times
\\
&\begin{split}
\times H(\frac{y}{\sqrt{1-\gamma}})&\left[
-(1-\gamma)^{\frac 32} \left(\frac{x_1-y_1}{\sqrt{\gamma}}H_1(\frac{x-y}{\sqrt{\gamma}})+\frac{x_2-y_2}{\sqrt{\gamma}}
H_2(\frac{x-y}{\sqrt{\gamma}})\right)
+\right.
\\
 & \left. +(\gamma(1-\gamma))^{\frac 12}\left(y_1H_1(\frac{x-y}{\sqrt{\gamma}})+y_2H_2(\frac{x-y}{\sqrt{\gamma}})\right)
\right ] dy
\end{split}
\end{split}
\end{equation}
Ya.G. Sinai and Dong Li in \cite{DS} show that operator acting on functions $\{H_1,H_2\}$ given by right-hand side of the equation \eqref{transformedeq} has $6$ unstable eigen values and $4$ neutral. In other words for any $T_{cr}>0$ there exist such constant $A(T_{cr})$ and $10$-parametric family of initial data $c_0(k)$ such that corresponding functions $g_p$ after renormalisation tends to some limit form $\tilde{g}(y,s)$,  while $p$ tends to infinity.  
Energy of the solution of the equation \eqref{eq: NSS}, constructed from functions $\tilde{g}(y,s)$ tends to infinity as $t\to T_{cr}$. In \cite{DS} there was obtained an asymptotics for energy growth
\begin{displaymath}
\mathcal{E}(t)=\int\limits_{\mathbb{R}^3}|v(k,t)|^2dk\simeq \frac{1}{(T_{cr}-t)^5}
\end{displaymath}
Proof of the existence theorem in \cite{DS} does not contain particular form of the initial data and corresponding solution. Thus numerical study of such solutions is of known interest.

\section{Computational scheme and main result.} 
Differentiation by time $t$ in the lefthand side of the equation \eqref{eq: NSS} represents traditional form of the ordinary differential equation for the vector-valued function $v(k,t)$ for any fixed $k$. Convolution-type operator in the righthand side has non-local nature i.e. value of solution $v(k,t)$ in particular point $(k,t)$ depends on values of the solution in each point $k \in \mathbb{R}^3$. Computational model is constructed in such domain $\Omega\subset \mathbb{R}^3$ that for all $0\leqslant t <T_{cr}$ values of function $v(k,t)$ outside of the  $\Omega$ can be neglected. For the discretisation of the arguments $t$ and $k$ we consider difference operator for $t$ and computation with discrete argument instead of integration in the righthandside of \eqref{eq: NSS}

Rewrite space integration term in the equation \eqref{eq: NSS} using incompressibility condition 
\begin{displaymath}
\int\limits_{\mathbb{R}^3}\langle k,v(k-l,s)\rangle P_k v(l,s)dl=
P_k\int\limits_{\mathbb{R}^3}\langle l,v(k-l,s)\rangle v(l,s)dl
\end{displaymath}
Expression in the righthandside can be treated in terms of convolutions of three-dimensional arrays.

Look for solutions of the equation \eqref{eq: NSS} with initial data located in the ball of radius $O(\sqrt{R})$ centered in the point $k_0=(0,0,R)$ and having in this support form of linear combination of rescaled Hermit functions. In more details we compute initial values  $v(k,0)=(v^1(k),v^2(k),v^3(k))$ in points of three-dimensional grid with given accuracy using expression
             
\begin{equation}
\label{eq: initial_conditions}
v^j(k)=\prod_{i=1}^{3}\sum\limits^{D}_{m=1}\lambda_{ijm}
He^{(m)}\left(\frac{3(k^i-k_0^i)}{\sqrt{R}}\right),\quad j=1,2,3
\end{equation}
Here $k^i$, $v^j$ correspond to components of the vectors $k$ and $v$ respectively, $He^{(m)}\left(\frac{3(k^i-k_0^i)}{\sqrt{R}}\right)$ is a Hermite function of order 
$m$ with argument $\frac{3(k^i-k_0^i)}{\sqrt{R}}$ and $\left\{\lambda_{ijm}\right\}$ corresponds  to the set of coefficients.
Available computational resources made it possible to consider grids of order $(\const\cdot10^2)\times(\const\cdot10^2)\times(\const\cdot10^3)$, where constants may vary from 
$1$ to  $2$. According to this accuracy of grid in each coordinate vary from $0.1$ to $0.2$ respectively.
For the case $R=5$ support of the solution $v(k,t)$, which we were able to observe lie in the domain $[-15,15]\times[-15,15]\times [1,200]$.

For calculation of three-dimensional convolutions we used Intel freeware library for non-commercial using \emph{Math Kernel
Library}.    

Since due to \cite{S1} for the solution $v(k,t)$ we have asymptotic $|k|\sqrt{t}\}$ one can estimate accuracy of the computational scheme. With space resolution near $0.2$ relative accuracy can be estimated with $10\%$.

\subsection*{Tornado--type solutions.}
{\bf Main result of this paper consists in presenting particular form of the initial conditions for which energy of the solution $\mathcal{E}(t)$ has polynomially fast growthrate while $t$ tends to $T_{cr}$}.       
 
Rough, due to our computational resources, estimate of the order of growth rate of the energy of the solution is far from the teoretical ${1}/{(T_{cr}-t)^5}$. Observed asymptotics for this order with $T_{cr}$ defined with respect to accuracy of our calculations are near $15$ -- $20$. Thus one may suppose that formation of blowing up solution, consturted from the functions $\tilde{g}(y,s)$ as it proposed in \cite{DS} is not unique. 

For the moment we have no examples of the solutions with growth rate $\sim 5$. For large set of variations of parameters $D$, $\lambda_{ijm}$ and form of initial condition solution decay in time.

However, founded behavior of numerical solution is stable for slightly perturbed (near $25\%$) discretisation step in time, space and form of initial data. Form of the blowing solution for such perturbed initial data differs very slightly, and time of the solution transforms with natural rescaling of \eqref{eq: NSS}.

\section{Results.}
Below we present as illustration computed results for the support located in the ball of radius $R=5$, centered at the point $\left(0,0,5\right)$. Maximal order of Hermite polynom is taken to be $D=3$, starting discretisation time step is $dt=0.001$, space resolution is $0.14907$ thus number of points of discretisation in X and Y coordinates are 
$nX=nY=160$, number of discretisation points on third coordinate is $nZ=1440$. We present square root of the integrated square of the absolute value of the solution ($mynorm=\sqrt{\mathcal{E}}$) in consequent time moments. From $t=0.000$ to $t=0.042$ energy monotonically decreased. 

\begin{tabular}[t]{ ll|ll|ll|ll}
t=0.000 & $\mynorm(t)$=1.00000 & t=0.042 & $\mynorm(t)$=0.3664 &t=0.045 & $\mynorm(t)$=19.972 \\
t=0.040 & $\mynorm(t)$=0.32191 & t=0.043 & $\mynorm(t)$=0.9500 &t=0.046 & $\mynorm(t)$=101.20 \\
t=0.041 & $\mynorm(t)$=0.31746 & t=0.044 & $\mynorm(t)$=4.1879 &t=0.047 & $\mynorm(t)$=698.57
\end{tabular}

In the second table we present results of the computations with the same initial data with different discretisation parameters. Space discretisation here is $0.20328$, time step $dt=0.0008$, number of points in X and Y coordinartes are $nX=nY=120$ and for third coordinate $nZ=1080$. 

\begin{tabular}[t]{ll|ll|ll}
 \dots      & \ldots                   & t=0.0424 & $\mynorm(t)$=0.3047     & t=0.0456 & $\mynorm(t)$=1.1492 \\
 t=0.0400   & $\mynorm(t)$=0.3226     & t=0.0432 & $\mynorm(t)$=0.2996     & t=0.0464 & $\mynorm(t)$=4.6551  \\
 t=0.0408   & $\mynorm(t)$=0.3165     & t=0.0440 & $\mynorm(t)$=0.3017     & t=0.0472 & $\mynorm(t)$=20.669  \\
 t=0.0416   & $\mynorm(t)$=0.3105     & t=0.0448 & $\mynorm(t)$=0.3982     & t=0.0480 & $\mynorm(t)$=1027.7  
\end{tabular}

\begin{figure}
\includegraphics[scale=0.8]{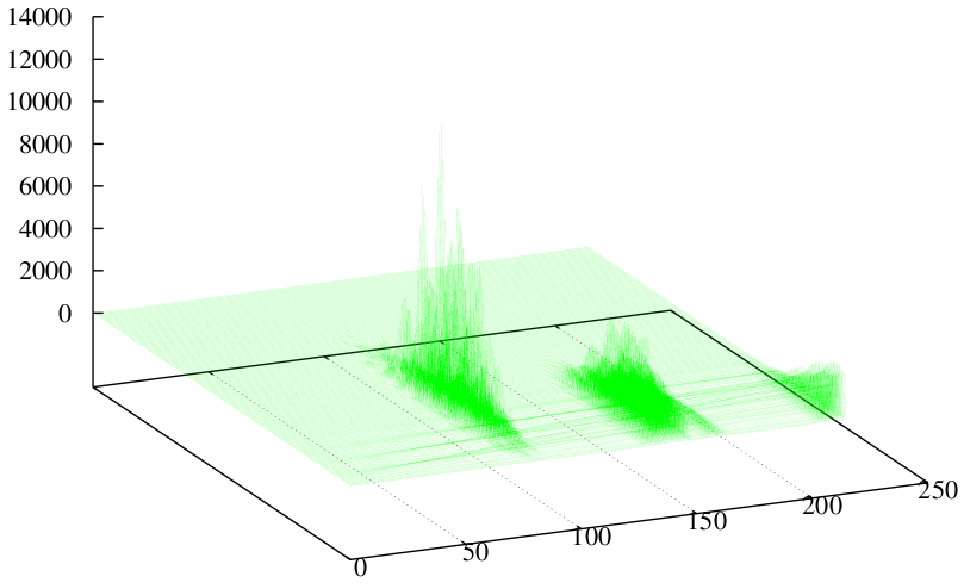}
\includegraphics[scale=0.8]{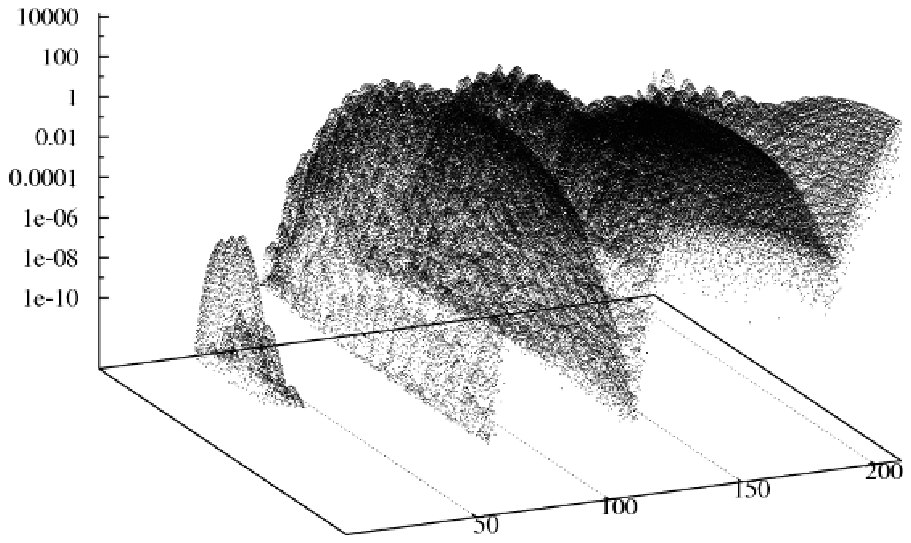}
\caption{Absolute Value in uniform and log scales}
\end{figure}

\begin{figure}
\includegraphics[scale=0.8]{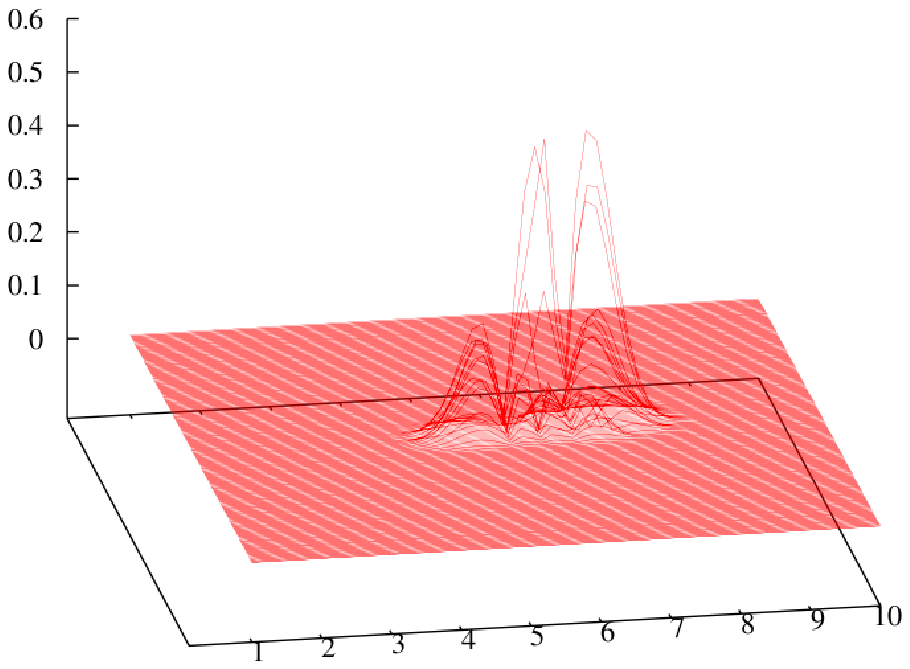}
\includegraphics[scale=0.8]{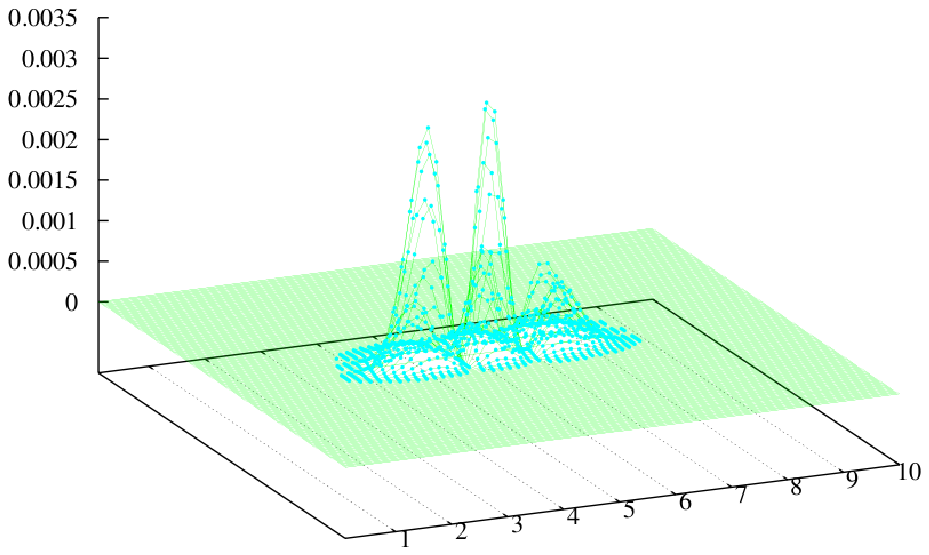}
\caption{Initial Value at t=0 and corresponding part of the solution at t=0.047}
\end{figure}

\begin{figure}
\includegraphics[scale=0.8]{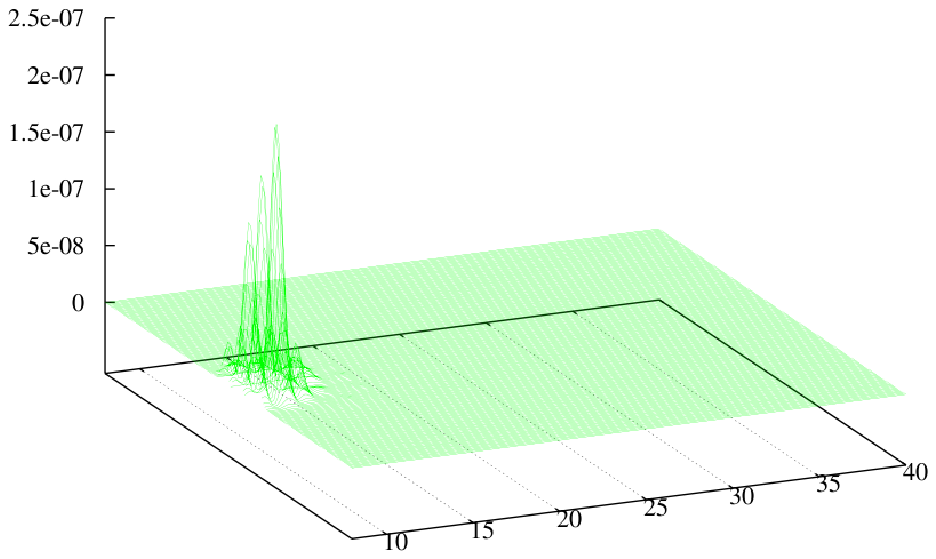}
\includegraphics[scale=0.8]{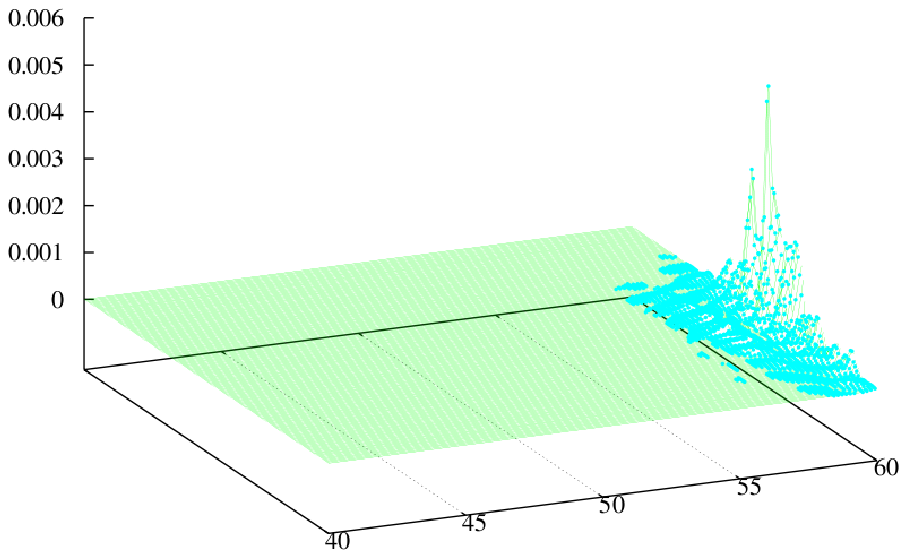}
\caption{Almost zeros for $z\in [20,60]$}

\includegraphics[scale=0.8]{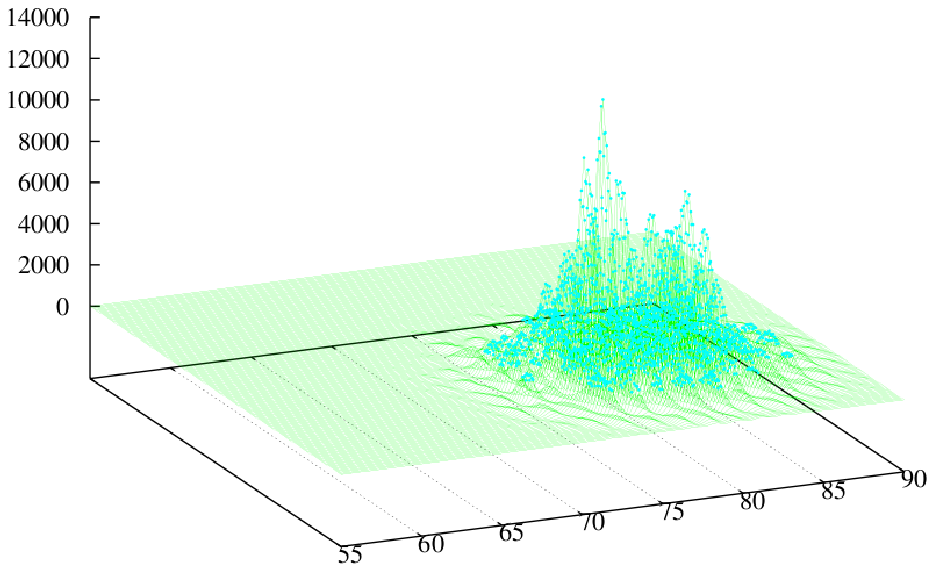}
\includegraphics[scale=0.8]{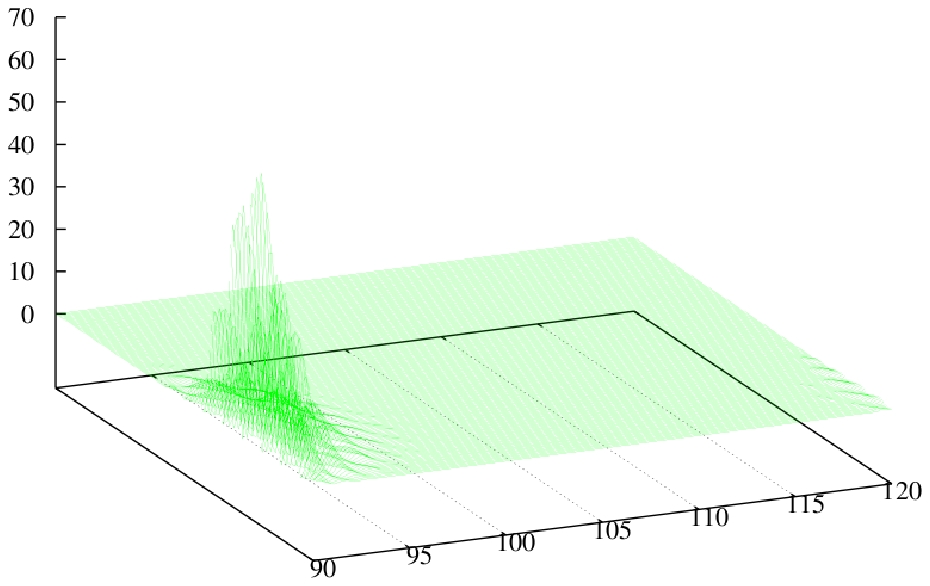}
\caption{First cloud. Main contribution.}

\includegraphics[scale=0.8]{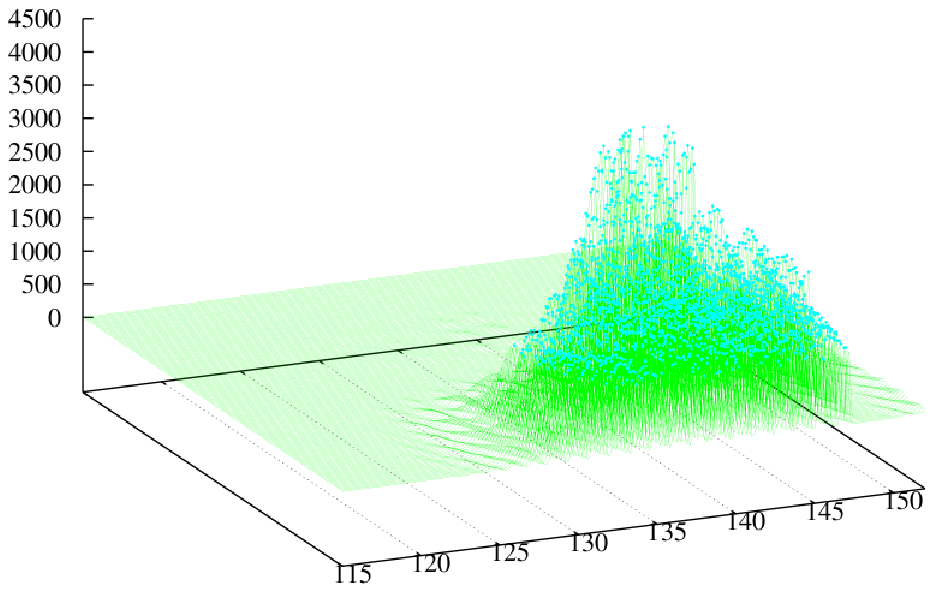}
\includegraphics[scale=0.8]{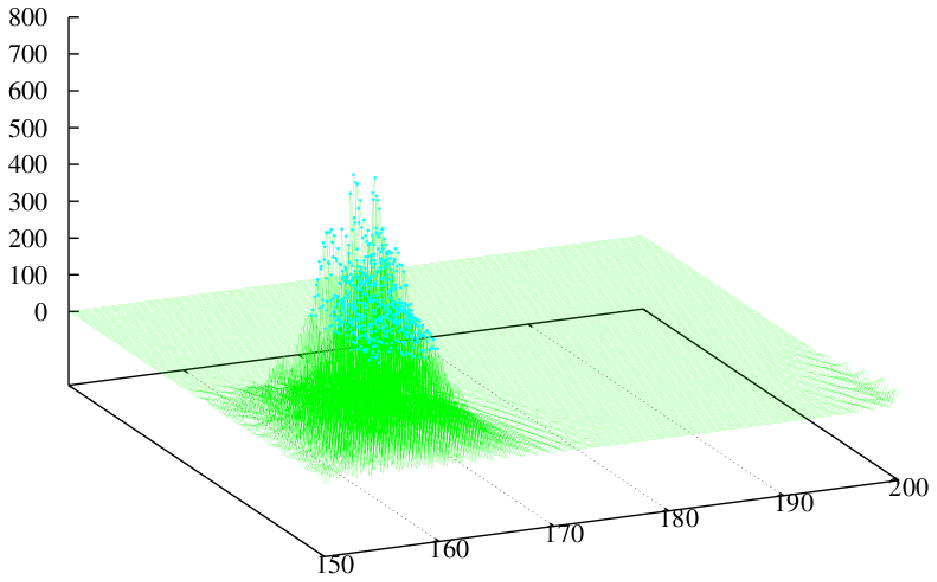}
\caption{Second and Third clouds.}
\end{figure}

\end{document}